\documentclass[10pt]{article}
\usepackage{amsmath}
\usepackage{amssymb}
\usepackage{mathrsfs}
\usepackage{amsfonts}
\usepackage{mathrsfs,amscd,amssymb,amsthm,amsmath,bm,graphicx}
\usepackage{graphicx}
\usepackage{cite}
\usepackage{color}
\usepackage{enumitem}
\usepackage{setspace}
\makeatletter

\renewcommand{\@seccntformat}[1]{{\csname the#1\endcsname}{\normalsize.}\hspace{.5em}}
\makeatother

\def \[{\begin{equation}}
\def \]{\end{equation}}

\newtheorem{thm}{Theorem}[section]

\newtheorem{lem}[thm]{Lemma}

\newtheorem{cor}[thm]{Corollary}

\begin{document}
\setlength{\baselineskip}{13pt}
\begin{center}{\Large \bf
The Laplacian spectrum, Kirchhoff index and complexity of the linear heptagonal networks
}

\vspace{4mm}

{\large Jia-Bao Liu$^{1}$, Jing Chen$^1$,Jing
Zhao$^1$, Shaohui Wang$^{2,*}$}\vspace{2mm}

{\small $^1$School of Mathematics and Physics, Anhui Jianzhu
University, Hefei 230601, P.R. China}
\\ {\small $^2$Department of Mathematics, Louisiana College, Pineville, LA 71359, USA}
\vspace{2mm}
\end{center}

\footnotetext{E-mail address: liujiabaoad@163.com,
chenjing102810@163.com, zhaojing94823@163.com, shaohuiwang@yahoo.com.}

\footnotetext{* Corresponding author.}

 {\noindent{\bf Abstract.}\ \ Let $H_n$ be the linear heptagonal networks with $2n$ heptagons. We study the structure properties and the eigenvalues of the linear heptagonal networks. According to the Laplacian polynomial of $H_n$, we utilize the decomposition theorem. Thus, the
 Laplacian spectrum of $H_n$ is created by eigenvalues
of a pair of matrices: $L_A$ and $L_S$ of order number $5n+1$ and $4n+1$, respectively. On the basis of the roots
and coefficients of their characteristic polynomials of $L_A$ and $L_S$, we not only get the explicit forms of Kirchhoff index, but also corresponding total complexity of $H_n$.

\vspace{2mm}

}

\section{Introduction}
\ \ \ \ With the discovery of carbon nanotubes, there has been a
practical interest in the study of concerning structures based on hexagon networks. Especially for a class of linear hexagonal networks, Yang \cite{MA,Mohar} and Huang et al. \cite{HYF,H.H} respectively explored Kirchhoff and another form multiplicative degree-Kirchhoff index of the linear hexagonal networks. Years later, Peng \cite{Y.J,J.L} et al. had studied the Kirchhoff index and related complexity for the linear phenylenes which are constructed by polyominos and hexagons. Inspired by \cite{Zhu,W.I}, the multiplicative version of degree-Kirchhoff index and its complexity of the generalized phenylenes we given. We refer the associated work \cite{C,A} and the notation we use in this paper follow \cite{J.H,A.D}.

Tips of carbon nanotubes are closed by pentagons. Heptagons, which generate negative curvatures, can also
be introduced with those pentagons \cite{H.Y,S.H}. Additionally, there are abundant practical applications of the heptagons in nature aspects, a remarkable example is that cacti are the most common plants with heptagons in natural structures. In our work, we concern a class of the linear heptagonal networks $H_n$ which inspired by these linear hexagonal networks, see Figure 1. Evidently, we see $\pi$ is an automorphism $\pi=(\bar{1})(\bar{2})\cdots(\bar{n})(1,\widehat{1})(2,2')\cdots(4n+1,(4n+1)')$
.

Before proceeding, we shall disgress to review some terminologies for the chemical graph theory. Since the topological indices are closely related to
the physical and chemical properties of the corresponding molecular graph, the calculation of various topological indices
is the core subject of chemical graph theory. At this point, we slightly to observe some properties of the Kirchhoff index. Klein and Randi\'{c} in \cite{D.J,J} defined a summation of all resistance distances between each pairs of vertices from $G$ as its Kirchhoff index $Kf(G)$, namely
\begin{eqnarray*}
Kf(G)=\sum_{i<j}r_{ij},
\end{eqnarray*}
where $r_{ij}$ is its resistance distance between two vertices $i$ and $j$.

Lately, one finds Kirchhoff index of a graph $G$ is closely related to the Laplacian eigenvalue~\cite{J.H,Gut,K.L}. That is
\begin{eqnarray*}
Kf(G)=n\sum_{k=2}^n\frac{1}{\mu_k},
\end{eqnarray*}
in which $0=\mu_1<\mu_2\leq\cdots \leq \mu_n$ are Laplacian eigenvalues of
$G$.

The organization of the paper is under the order. For section 2, some necessary notation of the block matrices are reviewed. Then, we apply automorphism and some lemmas into $H_n$. Lately, the Laplacian matrix can be consisted of blocks matrices, we can obtain the results of block matrices. Finally, we use transposition of matrices
to get $L_A$ and $L_S$. In section 3, we first study those explicit formulas of Kirchhoff index for $H_n$. Connecting these eigenvalues of $L_A$ and $L_S$. We can obtain the Laplacian spectrum of $H_n$. In the end, we have obtained complexity of $H_n$.
\begin{figure}[htbp]
\centering\includegraphics[width=16cm,height=6cm]{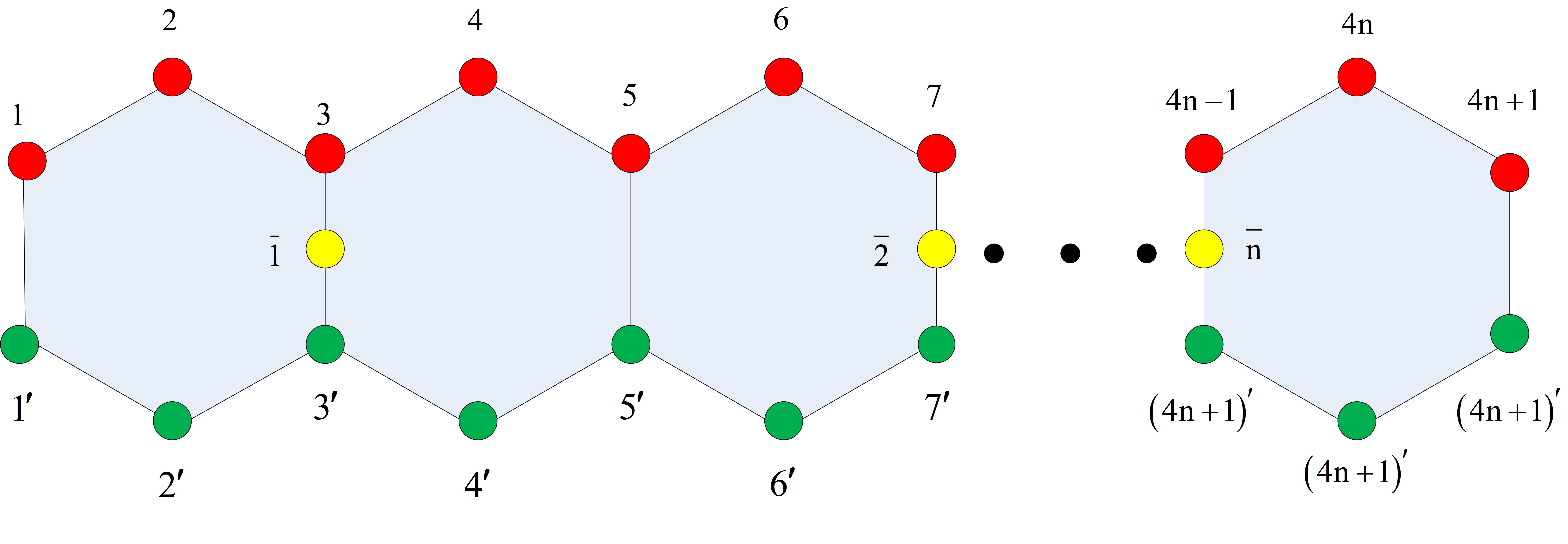}
\caption{The networks $H_n$.}
\end{figure}

\section{Preliminaries}
\ \ \ \ Before on embarking on the proofs of our main results, we disgress straightforward to observe a special method that are the basis of our
work. These equations that we describe are referred to~\cite{C,YY,P.A}.

For an $n$ by $n$ matrix $M$, denote by this submatrix of
$M$ the $M[s_1,\ldots,s_k]$, and created by the deletion of the
$s_1$-th, \ldots, $s_k$-th rows and corresponding columns. In the following part, $P_{M}(x)=\det (xI-M)$ is this characteristic polynomial of
$M$.

Label these vertex set of $H_n$ by
$V_0=\{\bar{1},\bar{2},\ldots,\bar{n}\}$, $V_1=\{1,2,\ldots,4n+1\}$
and $V_2=\{\widehat{1},2',\ldots,(4n+1)'\}$. So the Laplacian matrix of $H_n$
will be block matrices below~\cite{KDJ,LUK}.
\begin{equation}
L(H_n)=\left(
 \begin{array}{ccc}
  L_{V_0V_0}&
  L_{V_0V_1}&L_{V_0V_2}\\
  L_{V_1V_0}&L_{V_1V_1}& L_{V_1V_2}\\
 L_{V_2V_0}&L_{V_2V_1}& L_{V_2V_2}\\
  \end{array}
\right).
\end{equation}

Attention that, $L_{V_0V_1}=L_{V_0V_2},
L_{V_1V_0}=L_{V_2V_0}, L_{V_1V_1}=L_{V_2V_2}$
and $L_{V_1V_2}=L_{V_2V_1}$.

Given
\begin{equation*}
T=\left(
  \begin{array}{ccc}
I_{n}& 0&0\\
  0&\frac{1}{\sqrt{2}}I_{n}& \frac{1}{\sqrt{2}}I_{n}\\
 0& \frac{1}{\sqrt{2}}I_{n}& -\frac{1}{\sqrt{2}}I_{n}
  \end{array}
\right).
\end{equation*}

Consider the unitary transformation of $TL(H_n)T'$, then
\begin{equation}
TL(H_n)T'=\left(
  \begin{array}{cc}
  L_A& 0\\
    0& L_{S}
  \end{array}
\right),
\end{equation}
where $T'$ is the transposition of $T$. Obviously,

\begin{equation}
L_A=\left(
  \begin{array}{cc}
I_{n}& \sqrt{2}L_{V_0V_1}\\
  \sqrt{2}L_{V_1V_0}&L_{V_1V_1}+L_{V_1V_2}\\
  \end{array}
\right),~L_S=L_{V_1V_1}-L_{V_1V_2}.
\end{equation}

What's more, we at this point present those matrices introduced above in the
following. According to the structure of {Figure 1}, one obtains
$L_{V_0V_0}=2I_n$.
\begin{eqnarray*}
\L_{V_0V_1}&=&(l_{st})_{n\times(4n+1)}\\&=&
\left(
  \begin{array}{ccccccccccc}
    0& 0& -1& 0& 0& 0& 0 & \cdots & 0 & 0 & 0\\
    0& 0& 0& 0& 0& 0& -1 & \cdots & 0 & 0 &  0\\
    0& 0& 0& 0& 0& 0& 0 & \cdots & 0 & 0 &  0\\
    \vdots& \vdots& \vdots& \vdots& \vdots& \vdots& \vdots & \ddots & \vdots & \vdots &  \vdots\\
    0& 0& 0& 0& 0& 0& 0 & \cdots & 0 & 0 &  0\\
    0& 0& 0& 0& 0& 0& 0 & \cdots & -1 & 0 & 0 \\
  \end{array}
\right)_{n\times (4n+1)},
\end{eqnarray*}
where $l_{st}=-1$ and $l_{st}=0$, else $t=4s-1$.
\begin{eqnarray*}
 L_{V_1V_1}=
\left(
  \begin{array}{ccccccccccc}
    2 & -1&0 &0&0  &\cdots &0 &0 &0 &0 &0  \\
    -1& 2 & -1&0 &0&\cdots &0  &0 &0 &0 &0\\
    0 & -1 & 3 & -1 &0&\cdots &0 &0 &0 &0 &0\\
    0 &0 & -1&2 & -1&\cdots &0  &0 &0 &0 &0\\
    0 &0 & 0&-1 & 3&\cdots &0  &0 &0 &0 &0\\
 \vdots &\vdots&\vdots &\vdots& \vdots&\ddots &\vdots  &\vdots &\vdots &\vdots &\vdots \\
     0 & 0& 0& 0& 0&\cdots&3 &-1 &0 &0 &0 \\
    0 & 0& 0& 0& 0&\cdots&-1 &2 &-1 &0 &0 \\
     0 & 0& 0& 0& 0&\cdots&0 &-1&3 &-1 &0 \\
     0 & 0& 0& 0& 0&\cdots&0 &0&-1 &2 &-1 \\
     0 & 0& 0& 0& 0&\cdots&0 & 0&0 &-1 &2 \\
  \end{array}
\right)_{(4n+1)\times (4n+1)}.
\end{eqnarray*}

In addition, one gets
$$L_{V_1V_2}=(l_{st})_{(4n+1)\times(4n+1)}=diag(-1,
0, 0, 0, -1, \ldots, -1, 0, 0, 0,
-1),$$ where
$l_{4z+1,4z+1}=-1$ for $z\in\{0,1,2,\ldots,n\}$ and else $0$.

According to Equation {(2.3)}, one has
\begin{eqnarray*}
L_{V_1V_1}+
 L_{V_1V_2}&=&(l_{it})_{(4n+1)\times(4n+1)}\\
 &=&
\left(
  \begin{array}{ccccccccccccccc}
   1 & -1 &0 &0 & 0&0 &0 &\cdots &0 &0 &0 &0 &0 &0 &0 \\
    -1 & 2 & -1& 0 &0 &0&0&\cdots &0  &0 &0 &0 &0 &0 &0\\
    0 & -1 & 3& -1&0&0&0 &\cdots &0 &0 &0 &0 &0 &0 &0 \\
    0 &0 & -1 &2 &-1&0&0&\cdots &0  &0 &0&0&0&0&0\\
    0 &0 & 0&-1 &2 &-1&0&\cdots &0  &0 &0&0&0&0&0\\
    0 &0 & 0&0 &-1 &2&-1&\cdots &0  &0 &0&0&0&0&0\\
    0 &0 & 0&0 &0 &-1&3&\cdots &0  &0 &0&0&0&0&0\\
    \vdots&\vdots&\vdots &\vdots &\vdots &\vdots& \vdots&\ddots &\vdots   &\vdots&\vdots&\vdots&\vdots&\vdots&\vdots\\
     0 & 0 &0 &0 & 0&0 &0 &\cdots &3 &-1 &0 &0 &0 &0 &0 \\
     0 & 0 &0 &0 & 0&0 &0 &\cdots &-1 &2 &-1 &0 &0 &0 &0 \\
     0 & 0 &0 &0 & 0&0 &0 &\cdots &0 &-1&2 &-1 & 0&0 &0 \\
     0 & 0 &0 &0 & 0&0 &0 &\cdots &0 &0 &-1 &2 &-1 &0 &0 \\
     0 & 0 &0 &0 & 0&0 &0 &\cdots &0 &0 &0 &-1 &3 &-1 &0 \\
     0 & 0 &0 &0 & 0&0 &0 &\cdots &0 &0 &0 &0 &-1 &2 &-1 \\
     0 & 0 &0 &0 & 0&0 &0 &\cdots &0 &0 &0 &0 &0 &-1 &1 \\
  \end{array}
\right)_{(4n+1)\times (4n+1)},
\end{eqnarray*}
where $l_{11}=l_{4n+1,4n+1}=1$, $l_{ss}=2$ with
$s\in\{2,4,\ldots,4n\}$ and for $s\neq4z-1$,
$l_{4z-1,4z-1}=3$ for $z\in\{1,2,\ldots,n\}$. Also,
$l_{s,s+1}=l_{s+1,s}=-1$ for
$s\in\{1,2,3,\ldots,4n\}$ and $l_{st}=0$ wsth $|s-t|>1$.

And
\begin{eqnarray*}
L_{S}&=&(l_{st})_{(4n+1)\times(4n+1)}\\&=&
\left(
  \begin{array}{cccccccccc}
    3 & -1 &0 &0 & 0&\cdots &0 &0 &0 \\
    -1 & 2 & -1& 0 &0 &\cdots &0  &0 &0\\
    0 &-1 & 3 & -1&0 &\cdots &0 &0 &0 \\
    0 &0 & -1&2 &-1 &\cdots &0  &0 &0\\
    0 &0 & 0&-1 & 3 &\cdots &0  &0 &0\\
    \vdots &\vdots &\vdots &\vdots& \vdots&\ddots &\vdots  &\vdots &\vdots\\
     0 & 0& 0& 0& 0&\cdots&3 &-1 &0 \\
    0 &0 & 0& 0&0&\cdots &-1 &2 &-1  \\
     0 &0 &0 &0 & 0&\cdots&0 &-1 &3 \\
  \end{array}
\right)_{(4n+1)\times (4n+1)},
\end{eqnarray*}
where  $l_{ss}=3$ for
$s$ is odd and $l_{ss}=2$ for
$s$ is even
$l_{s,s+1}=l_{s+1,s}=-1$ with
$s\in\{1,2,3,\ldots,4n\}$ and $l_{st}=0$ with $|s-t|>1$.

Next, we list some lemmas as basic that may be used
throughout our main proof.
\begin{lem}~\cite{YW} Suppose that $L_A, L_S$  are mentioned before. We have
$$P_{L(H_n)}(x)=P_{L_A}(x)P_{L_{S}}(x).$$
\end{lem}
Based on the properties of these spanning trees and laplacian eigenvalues~\cite{PAL}, we have
the following lemma.
\begin{lem}~\cite{FR} Let $G$ be a connected graph with $|G|=n$ and $\tau(G)$ be the complexity of $G$. We have
$$\tau(G)=\frac{1}{n}\prod_{s=2}^{n}u_s.$$
\end{lem}

\begin{lem}

~\cite{ZF} Let $M_1, M_2, M_3$ and $M_4$ be the square matrices with $p\times p$ by $p\times q$,
$q\times p$ by $q\times q$ $q\times q$ by $q\times q$ and $q\times q$ by $q\times q$ matrices, respectively. We have
\begin{eqnarray*}
~det ~\left(%
\begin{array}{cc}
  M_1 & M_2 \\
  M_3 & M_4 \\
\end{array}%
\right)= det(M_1) \cdot det \big (M_4 - M_3M_1^{-1}
M_2\big),\\\nonumber
 \end{eqnarray*}
\end{lem}
where $M_1$ is invertible, and
$M_4-M_3M_1^{-1}M_2$ are denoted as Schur complements of $M_4$.
\section{Main results of $H_n$}
\ \ \ \ So far, our first main result provides the
explicit formula to the Kirchhoff index of $H_n$. Together with
their eigenvalues of $L_A$ and $L_S$, it is not hard
to use Lemma 2.1 and obtain Laplacian spectrum of $H_n$. Our second main result is to furnish the complete information of complexity of $H_n$, which
contributed by degree products and their eigenvalues. Suppose $0=\alpha_1< \alpha_2\leq \cdots\leq
\alpha_{5n+1}$ and $0<\beta_1\leq \beta_2\leq \cdots\leq
\beta_{4n+1}$ are respectively the roots of $P_{L_{A}}(x)=0$ and
$P_{L_{S}}(x)=0$. From (1.1), We obtain
\begin{lem}
Let $H_n$ be the linear heptagonal networks. We have
\begin{eqnarray}
Kf(H_n)&=&(9n+2)\cdot(\sum_{s=2}^{5n+1}\frac{1}{\alpha_s}+\sum_{t=1}^{4n+1}\frac{1}{\beta_t}).
\end{eqnarray}
\end{lem}

Obviously, we only need to compute the eigenvalues of the
$L_{A}$ and $L_{S}$. Thus, the formulas of
$\sum_{s=2}^{5n+1}\frac{1}{\alpha_s}$ and
$\sum_{t=1}^{4n+1}\frac{1}{\beta_t}$ are given as follows.

\begin{lem}
Assume that $0=\alpha_1< \alpha_2\leq \cdots\leq \alpha_{5n+1}$
are the eigenvalues of $L_{A}$. One has
$$\sum_{s=2}^{5n+1}\frac{1}{\alpha_s}=\frac{108n^3+66n^2+7n+4}{36n+8}.$$
\end{lem}
\noindent{\bf Proof.}~It is known that
$$P_{L_{A}}(x)=x^{5n+1}+a_1x^{5n}+\cdots+a_{5n-1}x^{2}+a_{5n}x=x(x^{5n}+a_1x^{5n-1}+\cdots+a_{5n-1}x+a_{5n})$$
with $a_{5n}\neq0 .$ Then, we find that $\frac{1}{\alpha_2},
\frac{1}{\alpha_3},\ldots, \frac{1}{\alpha_{5n+1}}$ satisfy the equation
$$a_{5n}\cdot x^{5n}+a_{5n-1}\cdot
x^{5n-1}+\cdots+a_{1}\cdot x+1=0.$$

By Vieta theorem, it gives
\begin{eqnarray*}
\sum_{s=2}^{5n+1}\frac{1}{\alpha_s}=\frac{(-1)^{5n-1}\cdot
a_{5n-1}}{(-1)^{5n}\cdot a_{5n}}.
\end{eqnarray*}
\noindent{\bf Proposition 1.}
 $(-1)^{5n}\cdot
a_{5n}=(9n+2)\cdot2^{n-1}.$\\
\noindent{\bf Proof.}
Because $(-1)^{5n}a_{5n}$ is the
summation of these principal minors from $L_{A}$, which have
$5n$ rows and columns, we get
\begin{eqnarray}
(-1)^{5n}\cdot a_{5n}=\sum_{s=1}^{5n+1}\det L_{A}[s].
\end{eqnarray}

In view of equation (2.5), we know that $s$ can be chosen from the
matrix $2I_n$ or
$L_{V_1V_1}+L_{V_1V_2}.$ Hence, we consider the cases below.

{\bf{Case 1.}} For the case of $1\leq s\leq n.$ In the other words, $L_{A}[s]$ is created by dropping the $s$
row of matrices $2I_n$ and
$\sqrt{2}L_{V_0V_1}$, and also the related columns of
 $2I_n$ and $\sqrt{2}L_{V_1V_0}.$ Let $2I_{n-1},~
R_{(n-1)\times(4n+1)},~R_{(n-1)\times(4n+1)}'$ and
$U_{(4n+1)\times(4n+1)}$ be these blocks.
By using Lemma 2.3 and fundamental operations, one knows
\begin{eqnarray*}
\det L_{A}[s]=\left|
  \begin{array}{cc}
    2I_{n-1} & R_{(n-1)\times(4n+1)} \\
    R_{(n-1)\times(4n+1)}' & U_{(4n+1)\times(4n+1)} \\
  \end{array}
\right|=\left|
  \begin{array}{cc}
    2I_{n-1} & 0 \\
    0 & M_1 \\
  \end{array}
\right|,
\end{eqnarray*}
where
$M_1=U_{(4n+1)\times(4n+1)}-R_{(n-1)\times(4n+1)}'R_{(n-1)\times(4n+1)}=(m_{pq}),$
of which $ m_{11}=m_{4n+1,4n+1}=1, m_{pp}=2$ for
$p\in\{2,3,\ldots,4n-1,4n\},$
 but $p\neq 4s-1,$ $m_{4s-1,4s-1}=3.$ Also,
$m_{p,p+1}=m_{p+1,p}=-1$ for
$p\in\{1,2,\ldots,4n+1\}$ and $m_{pq}=0$ for $|p-q|>1$. By routine calculations, we have
\begin{eqnarray}
\det L_{A}[s]=\det M_1=(2)^{n-1}.
\end{eqnarray}

{\bf{Case 2.}} For the case of $n+1\leq s\leq
5n+1.$ We say $L_{A}[s]$ is deduced by deleting the
$s$ row of the $L_{V_1V_1}+L_{V_1V_2}$ and
$\sqrt{2}L_{V_1V_0}$, also the corresponding column of
the $L_{V_1V_1}+L_{V_1V_2}$ and
$\sqrt{2}L_{V_0V_1}.$ Denote by $2I_{n},~
R_{n\times4n},~R_{n\times4n}'$ and $U_{4n\times4n}$ desired
blocks. By Lemma 2.3 and fundamental operations, we obtain that
\begin{eqnarray*}
\det L_{A}[s]=\left|
  \begin{array}{cc}
    2I_{n} & R_{n\times4n} \\
    R_{n\times4n}' & U_{4n\times4n} \\
  \end{array}
\right|=\left|
  \begin{array}{cc}
    2I_{n} & 0 \\
    0 & M_2 \\
  \end{array}
\right|,
\end{eqnarray*}
where $M_2=U_{4n\times4n}-R_{n\times4n}'R_{n\times4n}.$ Set
$M_{2'}=(m_{pq})$ with order $4n+1$, of which  $m_{pp}=2$ for
$p\in\{2,3,4,5,6,7,\ldots,4n-1,4n\}$. Also,
$m_{11}=m_{4n+1,4n+1}=1$,
$m_{p,p+1}=m_{p+1,p}=-1$ for
$p\in\{1,2,3,\ldots,4n-1,4n\}$ and $m_{pq}=0$ with $|p-q|>1$. By routine calculation, we have
\begin{eqnarray}
\det L_{A}[s]=2^n,
\end{eqnarray}

Combining Equations {(3.6)}--{(3.7)}, we have
\begin{eqnarray*}
(-1)^{5n}\cdot a_{5n}&=&\sum_{s=1}^{5n+1}\det
L_{A}[s]=\sum_{s=1}^{n}\det
L_{A}[s]+\sum_{s=n+1}^{5n+1}\det
L_{A}[s]\\
&=&n\cdot2^{n-1}+(4n+1)\cdot
2^n\\
&=&(9n+2)\cdot2^{n-1}.
\end{eqnarray*}

So the proof of the {Proposition
1} is finished.

Next, we will focus on recursive calculation for
$(-1)^{5n-1}a_{5n-1}$.\\

\noindent{\bf Proposition 2.}
 $(-1)^{5n-1}a_{5
n-1}=(108n^3+66n^2+7n+4)\cdot2^{n-3}.$\\
\noindent{\bf Proof.}
We know $(-1)^{5n-1}a_{5n-1}$ is
the summation of their principal minors of $L_{A}$ giving
$5n-1$ rows and columns, we can gain
\begin{eqnarray}
(-1)^{5n-1}\cdot a_{5n-1}=\sum_{s=1}^{5n+1}\det
L_{A}[s,t].
\end{eqnarray}

For the result of {Proposition 1}, we will consider the following cases:

{\bf{Case 1.}} For the case of $1\leq s< t\leq n.$
We see $L_{A}[ s,t]$ is obtained by dropping
$s$ and $t$ rows of the identity matrix $2I_n$ and
$\sqrt{2}L_{V_0V_1}$, also the corresponding columns of
the $2I_n$ and $\sqrt{2}L_{V_1V_0}.$ Denote $2I_{n-2},~
R_{(n-2)\times(4n+1)},~R_{(n-2)\times(4n+1)}'$ and
$U_{(4n+1)\times(4n+1)}$ by desired blocks.
Based on Lemma 2.3 and elementary operations, we
have
\begin{eqnarray*}
\det L_{A}[s,t]=\left|
  \begin{array}{cc}
    2I_{n-2} & R_{(n-2)\times(4n+1)} \\
    R_{(n-2)\times(4n+1)}' & U_{(4n+1)\times(4n+1)} \\
  \end{array}
\right|=\left|
  \begin{array}{cc}
    2I_{n-2} & 0 \\
    0 & M_3 \\
  \end{array}
\right|,
\end{eqnarray*}
where
$M_3=U_{(4n+1)\times(4n+1)}-R_{(n-2)\times(4n+1)}'R_{(n-2)\times(4n+1)}=(m_{pq}),$
of which $ m_{11}=m_{4n+1,4n+1}=1, $ $m_{pp}=2$ for
$p\in\{2,3,4,5,6,7,\ldots,4n-1,4n\}$ and $p\neq 4s-1$, also $p\neq 4t-1$,
$m_{4s-1,4s-1}=m_{4t-1,4t-1}=3.$ Also,
$m_{p,p+1}=m_{p+1,p}=-1$ with
$p\in\{1,2,3,\ldots,4n-1,4n\}$ with $m_{pq}=0$ with $|p-q|>1$. By
basic calculation, we have
\begin{eqnarray}
\det L_{A}[s,t]=\det
M_3=[4\cdot(t-s)+2]\cdot2^{n-2}.
\end{eqnarray}

{\bf{Case 2.}} For the case of $n+1\leq s<
t\leq 5n+1.$ We say $L_{A}[s,t]$ is created by
dropping $(s-n)$-th and $(t-n)$-th rows of the
$L_{V_1V_1}+L_{V_1V_2}$ and
$\sqrt{2}L_{V_1V_0}$, also the corresponding columns of
the $L_{V_1V_1}+L_{V_1V_2}$ and
$\sqrt{2}L_{V_0V_1}.$ Denote $2I_{n},~
R_{n\times(4n-1)},~R_{n\times(4n-1)}'$ and
$U_{(4n-1)\times(4n-1)}$ by desired blocks.
In terms of Lemma 2.3 and fundamental operations, we obtain
\begin{eqnarray*}
\det L_{A}[s,t]=\left|
  \begin{array}{cc}
    2I_{n} & R_{n\times(4n-1)} \\
    R_{n\times(4n-1)}' & U_{(4n-1)\times(4n-1)} \\
  \end{array}
\right|=\left|
  \begin{array}{cc}
    2I_{n} & 0 \\
    0 & M_4 \\
  \end{array}
\right|,
\end{eqnarray*}
where
$M_4=U_{(4n-1)\times(4n-1)}-R_{n\times(4n-1)}'R_{n\times(4n-1)}=M_{2'}[s-n,t-n]$. By recursive calculations, we obtain
\begin{eqnarray}
\det L_{A}[s,t]=
\begin{cases}
(t-1-n)2^{n-2},       & if~ s=n+1,~ t=n+z,~2\leq z \leq 4n; \\
(5n-s+1)2^{n},       & if~t=5n+1,~s=n+z,~2\leq z \leq 4n; \\
(t-s)2^{n},       & if~n+2\leq s<t\leq 5n.\\
\end{cases}
\end{eqnarray}

{\bf{Case 3.}} We consider the case of $1\leq s\leq n$,
$n+1\leq t\leq 5n+1$. In other words, we have $L_{A}[s,t]$ is created
by dropping $s$ rows of the identity matrix $2I_n$ and
$\sqrt{2}L_{V_0V_1}$,  $(t-n)$-th rows of the
$\sqrt{2}L_{V_1V_0}$ and
$L_{V_1V_1}+L_{V_1V_2}$, respectively. What is more, the corresponding columns of the $2I_n$,
$\sqrt{2}L_{V_1V_0}$,
$L_{V_1V_1}+L_{V_1V_2}$ and
$\sqrt{2}L_{V_0V_1}.$ Denote $2I_{n-1},~
R_{(n-1)\times4n},~R_{(n-1)\times4n}'$ and $U_{4n\times4n}$ by
desired blocks.By the methods of Lemma 2.3 and
fundamental operations, we have
\begin{eqnarray*}
\det L_{A}[s,t]=\left|
  \begin{array}{cc}
    2I_{n-1} & R_{(n-1)\times4n} \\
    R_{(n-1)\times4n}' & U_{4n\times4n} \\
  \end{array}
\right|=\left|
  \begin{array}{cc}
    2I_{n} & 0 \\
    0 & M_5 \\
  \end{array}
\right|,
\end{eqnarray*}
where
$M_5=U_{4n\times4n}-R_{(n-1)\times4n}'R_{(n-1)\times4n}=M_{1}[t-n]$.
By an explicit calculation, we get

\begin{eqnarray}
\det L_{A}[s,t]=
(|t-4s+1-n|+1)\cdot2^{n-1};
\end{eqnarray}

Considering equations {(3.9)}{--}{(3.11)}, one acquires
\begin{eqnarray*}
(-1)^{5n-1}a_{5n-1}&=&\sum_{1\leq s\leq
t}^{5n+1}\det L_{A}[s,t]\\
&=&\sum_{1\leq s\leq t}^{n}\det\ L_{A}[s,t]+ \sum_{n+1\leq
s\leq t}^{5n+1}\det L_{A}[s,t] +\sum_{1\leq s\leq n,
n+1\leq t\leq 5n+1}^{5n+1}\det\ L_{A}[s,t]\\
&=&\frac{4n^3+6n^2-7n}{6}\cdot2^{n-2}+\frac{256n^3+192n^2+32n}{6}\cdot2^{n-2}\\
&&+\frac{64n^3-4n+12}{6}\cdot2^{n-2}\\
&=&(108n^3+66n^2+7n+4)\cdot2^{n-3}.
\end{eqnarray*}

\begin{cor}
Suppose that $0=\alpha_1< \alpha_2\leq \cdots\leq \alpha_{4n+1}$ are the
roots of $P_{L_{A}}(x)=0$. Then
\begin{eqnarray*}
\sum_{s=2}^n\frac{1}{\alpha_s}&=&\frac{108n^3+66n^2+7n+4}{36n+8}.
\end{eqnarray*}
\end{cor}

Let
$P_{L_{s}}(x)=\det(xI-L_{s})=x^{4n+1}+c_1x^{4n}+\cdots+c_{4n}x+c_{4n+1}.$
By using the relation between the roots and coefficients of
$P_{L_{s}}(x)$, it is not hard to see
\begin{eqnarray*}
\sum_{t=1}^{4n+1}\frac{1}{\beta_z}=\frac{(-1)^{4n}c_{4n}}{(-1)^{4n+1}c_{4n+1}}=\frac{(-1)^{4n}c_{4n}}{\det L_{s}}.
\end{eqnarray*}

In the other way, we will concern $s$-th ordered principal submatrix, $M_s$, which is found by giving
first $s$ rows and columns of $ L_{s},
s=1,2,\cdots,n-1$. Let $m_s=\det M_s$.\\
\noindent{\bf Fact 1.} For $1\leq s\leq 4n$,\\
\begin{eqnarray*}
m_s=
\begin{cases}
\frac{1+\sqrt{3}}{2}\Big(\frac{\sqrt{2}+\sqrt{6}}{2}\Big)^{s}+\frac{1-\sqrt{3}}{2}\Big(\frac{\sqrt{2}-\sqrt{6}}{2}\Big)^{s}, & if~ s~  is ~ even;\\
\frac{\sqrt{6}+3\sqrt{2}}{4}\Big(\frac{\sqrt{2}+\sqrt{6}}{2}\Big)^{s}+\frac{-\sqrt{6}+3\sqrt{2}}{4}\Big(\frac{\sqrt{2}-\sqrt{6}}{2}\Big)^{s},  & if~ s~  is~  odd.\\
\end{cases}
\end{eqnarray*}
\noindent{\bf Proof.} We will compute and gain that $m_1=3,
m_2=5, m_3=12, m_4=19.$ For $3\leq s\leq 4n,$ we can extend $\det M_s $ with its last row as follows\\
\begin{eqnarray}
m_s=
\begin{cases}
2m_{s-1}-m_{s-2},& if~s~is~even;\\
3m_{s-1}-m_{s-2},& if~s~is~odd.\\
\end{cases}
\end{eqnarray}
for $1\leq s\leq 2n,$ let $e_s=m_{2s}$ and for $1\leq s\leq 2n-1,$let$ f_s=m_{2s+1}$. Then we have, together with (3.12), $e_s=5$, $f_s=12, $
for $s\geq2,$\\
\begin{eqnarray}
\begin{cases}
 e_s=2f_{s-1}-e_{s-1} ,& if~s~is~even;\\
f_s=3e_s-f_{s-1},& if~s~is~odd.\\
\end{cases}
\end{eqnarray}\\

From the first equation in (3.13), one has $f_{s-1}=\frac{e_s}{2}+\frac{e_{s-1}}{2},$ therefore $f_s=\frac{e_{s+1}}{2}+\frac{e_s}{2}.$
Replacing $f_{s-1}$ and $f_s$ into the second equation in (3.13) yield $e_{s+1}=4e_s-e_{s-1}$, $s\geq2$. Similarly, we obtain $f_{s+1}=4f_s-f_{s-1},$ $s\geq2$.
Thus, $\{m_s\}$ meets the recurrence relation\\
 \begin{eqnarray}
 m_1=3,m_2=5,m_3=12,m_4=19,m_s=4e_{s-2}-e_{s-4}.~(s\geq5)
 \end{eqnarray}

 So we see that the characteristic equation of (3.14) is $x^4=4x^2-1,$ whose roots are
 $x_1=\frac{\sqrt{2}+\sqrt{6}}{2},$ $x_2=-\frac{\sqrt{2}+\sqrt{6}}{2},$ $x_3=\frac{\sqrt{2}-\sqrt{6}}{2},$ $x_4=-\frac{\sqrt{2}+\sqrt{6}}{2}.$\\
Assume the general solution of (3.14) is\\
 \begin{eqnarray}
m_s=\Big(\frac{\sqrt{2}+\sqrt{6}}{2}\Big)^{s}y_1+\Big(-\frac{\sqrt{2}+\sqrt{6}}{2}\Big)^{s}y_2+
\Big(\frac{\sqrt{2}-\sqrt{6}}{2}\Big)^{s}y_3+\Big(-\frac{\sqrt{2}-\sqrt{6}}{2}\Big)^{s}y_4,
 \end{eqnarray}
combining the initial condition in (3.14), the system of equations are given\\
\begin{eqnarray*}
\begin{cases}
 \Big(\frac{\sqrt{2}+\sqrt{6}}{2}\Big)y_1+\Big(-\frac{\sqrt{2}+\sqrt{6}}{2}\Big)y_2+
 \Big(\frac{\sqrt{2}-\sqrt{6}}{2}\Big)y_3+\Big(-\frac{\sqrt{2}-\sqrt{6}}{2}\Big)y_4=3;\\
 \Big(\frac{\sqrt{2}+\sqrt{6}}{2}\Big)^{2}y_1+\Big(-\frac{\sqrt{2}+\sqrt{6}}{2}\Big)^{2}y_2+
 \Big(\frac{\sqrt{2}-\sqrt{6}}{2}\Big)^{2}y_3+\Big(-\frac{\sqrt{2}-\sqrt{6}}{2}\Big)^{2}y_4=5;\\
\Big(\frac{\sqrt{2}+\sqrt{6}}{2}\Big)^{3}y_1+\Big(-\frac{\sqrt{2}+\sqrt{6}}{2}\Big)^{3}y_2+
 \Big(\frac{\sqrt{2}-\sqrt{6}}{2}\Big)^{3}y_3+\Big(-\frac{\sqrt{2}-\sqrt{6}}{2}\Big)^{3}y_4=12;\\  \Big(\frac{\sqrt{2}+\sqrt{6}}{2}\Big)^{4}y_1+\Big(-\frac{\sqrt{2}+\sqrt{6}}{2}\Big)^{4}y_2+
 \Big(\frac{\sqrt{2}-\sqrt{6}}{2}\Big)^{4}y_3+\Big(-\frac{\sqrt{2}-\sqrt{6}}{2}\Big)^{4}y_4=19.\\
\end{cases}
\end{eqnarray*}\\

The unique solution of this system of equation is
$y_1=\frac{2+\sqrt{6}+2\sqrt{3}+3\sqrt{2}}{8},$ $y_2=\frac{2-\sqrt{6}+2\sqrt{3}-3\sqrt{2}}{8}$,
$y_3=\frac{2-\sqrt{6}-2\sqrt{3}+3\sqrt{2}}{8}$, $y_4=\frac{2+\sqrt{6}-2\sqrt{3}-3\sqrt{2}}{8}$.
Therefore Fact 1 deduces by (3.15).\\

By using the expansion formula, one can deduce $\det L_s $ with its last row below
$\det L_S =3M_{4n}-M_{4n-1}$.\\
Combining with Fact 1, it is immediately to see that \vspace{4mm}\\
\noindent{\bf Fact 2.} $\det L_s =\frac{(3+2\sqrt{3})(\sqrt{2}+\sqrt{6})^{4n}+(3-2\sqrt{3})(\sqrt{2}-\sqrt{6})^{4n}}{2^{4n+1}}$.\\
Now we are ready to determine $b_{4n}$. We can denote this diagonal of $\det L_s $ by $l_{ss},s=1,2,...,2n+1.$ \vspace{4mm}\\
\noindent{\bf Fact 3.} $b_{4n} =\frac{[120n+60\sqrt{3}n+6+5\sqrt{3}+(6+7\sqrt{3})(2-\sqrt{3})^{4n}](\sqrt{2}+\sqrt{6})^{4n}+
[120n-60\sqrt{3}n+6-5\sqrt{3}+(6-7\sqrt{3})(2+\sqrt{3})^{4n}](\sqrt{2}-\sqrt{6})^{4n}}{24\cdot2^{4n}}$.\\
\noindent{\bf Proof.}    Since $b_{4n}(=(-1)^{4n}b_{4n})$ is the summation of these principal minors from $L_s$ with size 4n by 4n, we have\\
\begin{eqnarray}
b_{4n}&=&\sum_{s=1}^{4n+1}\det
L_{S}[s]=\sum_{s=1}^{4n+1}\det\left(
  \begin{array}{cc}
    M_{s-1} & 0 \\
    0 & B_{1}\\
  \end{array}
\right)=\sum_{s=1}^{4n+1}\det M_{s-1}\det B_{1} ,
\end{eqnarray}\\
where
\begin{eqnarray*}
 B_{1}&=&
\left(
  \begin{array}{cccccc}
    l_{s+1,s+1} &-1 &\cdots & 0 \\
    \vdots & \ddots & \vdots & \vdots \\
    0 & \cdots & l_{4n,4n} &-1 \\
    0 & \cdots & -1&l_{4n+1,4n+1} \\
  \end{array}
\right).
\end{eqnarray*}\\

Because permutating rows and related columns in square matrices keep its determinant unchanged.
Combining with $L_{S}$, we will get $\det B_{1}=\det M_{4n+1-s}$. Therefore,\\
\begin{eqnarray}
b_{4n}&=&\sum_{s=1}^{4n+1}m_{s-1}m_{4n+1-s}=2m_{4n}+\sum_{z=1}^{2n-1}m_{2z}m_{4n-2z}+\sum_{z=1}^{2n}m_{2l-1}m_{4n-2l+1}.~(put:m_0=1)
\end{eqnarray}\\

By Fact 1, we have
$2m_{4n}=(1+\sqrt{3})\Big(\frac{\sqrt{2}+\sqrt{6}}{2}\Big)^{4n}+(1-\sqrt{3})\Big(\frac{\sqrt{2}-\sqrt{6}}{2}\Big)^{4n}$
and\\
\begin{eqnarray*}
\sum_{z=1}^{2n-1}m_{2z}m_{4n-2z}
&=&\sum_{z=1}^{2n-1}\Big[\frac{1+\sqrt{3}}{2}\Big(\frac{\sqrt{2}+\sqrt{6}}{2}\Big)^{2z}+
\frac{1-\sqrt{3}}{2}\Big(\frac{\sqrt{2}-\sqrt{6}}{2}\Big)^{2z}\Big]\\
&&\times \Big[\frac{1+\sqrt{3}}{2}\Big(\frac{\sqrt{2}+\sqrt{6}}{2}\Big)^{4n-2z}+
\frac{1-\sqrt{3}}{2}\Big(\frac{\sqrt{2}-\sqrt{6}}{2}\Big)^{4n-2z}\Big]\\
&=&\frac{[48n+24\sqrt{3}n-18-16\sqrt{3}+(6+4\sqrt{3})(2-\sqrt{3})^{4n}](\sqrt{2}+\sqrt{6})^{4n}}{24\cdot2^{4n}}\\&&+
\frac{[48n-24\sqrt{3}n-18+16\sqrt{3}+(6-4\sqrt{3})(2+\sqrt{3})^{4n}](\sqrt{2}-\sqrt{6})^{4n}}{24\cdot2^{4n}}.
\end{eqnarray*}
We can use a similar routine to obtain that
\begin{eqnarray*}
\sum_{l=1}^{2n}m_{2l-1}m_{4n-2l+1}
&=&\sum_{l=1}^{2n}\Big[\frac{\sqrt{6}+3\sqrt{2}}{4}\Big(\frac{\sqrt{2}+\sqrt{6}}{2}\Big)^{2l-1}+
\frac{-\sqrt{6}+\sqrt{3}}{4}\Big(\frac{\sqrt{2}-\sqrt{6}}{2}\Big)^{2l-1}\Big]\\
&&\times \Big[\frac{\sqrt{6}+3\sqrt{2}}{4}\Big(\frac{\sqrt{2}+\sqrt{6}}{2}\Big)^{4n-2l+1}+
\frac{-\sqrt{6}+3\sqrt{2}}{4}\Big(\frac{\sqrt{2}-\sqrt{6}}{2}\Big)^{4n-2l+1}\Big]\\
&=&\frac{[72n+36\sqrt{3}n-3\sqrt{3}+(9+6\sqrt{3})(2-\sqrt{3})^{4n+1}](\sqrt{2}+\sqrt{6})^{4n}}{24\cdot2^{4n}}\\&&+
\frac{[72n-36\sqrt{3}n+3\sqrt{3}+(9-6\sqrt{3})(2+\sqrt{3})^{4n+1}](\sqrt{2}-\sqrt{6})^{4n}}{24\cdot2^{4n}}.
\end{eqnarray*}\\
\begin{lem}

For graph $H_n,$\\
\begin{eqnarray*}
Kf(H_n)&=&(20n+2)\cdot\Big(\frac{108n^3+66n^2+7n+4}{36n+8}+\frac{(-1)^{4n}\cdot
b_{4n}}{\det L_{s}}\Big).
\end{eqnarray*}
\end{lem}
The Kirchhoff indices of linear heptagonal networks from
 $H_{1}$ to $H_{50}$ are denoted in Table 1 below.

\begin{lem}
For networks $H_n,$\\
\begin{eqnarray*}
\tau(L_n)&=&\frac{(3+2\sqrt{3})(\sqrt{2}+\sqrt{6})^{4n}+(3-2\sqrt{3})(\sqrt{2}-\sqrt{6})^{4n}}{2^{3n+2}}.
\end{eqnarray*}
{\bf Proof.} By Lemma 2.2, we have $
\tau(G)=\frac{1}{n}\prod_{s=2}^{n}u_s.
$ Note that
\begin{eqnarray*}
%
\prod_{s=2}^{5n+1}\alpha_s&=&(-1)^{5n}a_{5n}=(9n+2)(2)^{n-1},\\
%
\prod_{t=1}^{4n+1}\beta_t&=& \det
L_{s}=\frac{(3+2\sqrt{3})(\sqrt{2}+\sqrt{6})^{4n}+(3-2\sqrt{3})(\sqrt{2}-\sqrt{6})^{4n}}{2^{4n+1}}.
\end{eqnarray*}
\end{lem}

Then
\begin{eqnarray*}
\tau(H_n)&=&\frac{(3+2\sqrt{3})(\sqrt{2}+\sqrt{6})^{4n}+(3-2\sqrt{3})(\sqrt{2}-\sqrt{6})^{4n}}{2^{3n+2}}.
\end{eqnarray*}
Next, complexity in linear heptagonal networks from $H_{1}$ to $H_{12}$ are given in Table $2$.

\begin{table}[htbp]
\setlength{\abovecaptionskip}{0.05cm}
 \centering \vspace{.3cm}
\caption{Kirchhoff indices of linear heptagonal networks from $H_{1}$ to $H_{50}$.}
\begin{tabular}{cccccccccc}
  \hline
  $H_n$ & $Kf(H_n)$ & $H_n$ & $Kf(H_n)$ & $H_n$ & $Kf(H_n)$ & $H_n$ & $Kf(H_n)$ & $H_n$ & $Kf(H_n)$\\
  \hline
  $H_{1}$ & $79.25$ & $H_{11}$ & $41172.05$ & $H_{21}$ & $268961.03$ & $H_{31}$ & $845446.16$ & $H_{41}$ & $1932627.45$ \\
  $H_{2}$ & $404.17$ & $H_{12}$ & $52876.62$ & $H_{22}$ & $308245.21$ & $H_{32}$ & $928509.96$ & $H_{42}$ & $2075670.87$\\
  $H_{3}$ & $1138.09$ & $H_{13}$ & $66610.15$ & $H_{23}$ & $351178.36$ & $H_{33}$ & $1016842.73$ & $H_{43}$ & $2225603.25$\\
  $H_{4}$ & $2442.97$ & $H_{14}$ & $82534.64$ & $H_{24}$ & $397922.47$ & $H_{34}$ & $1110606.45$ & $H_{44}$ & $2382586.59$ \\
  $H_{5}$ & $4480.81$ & $H_{15}$ & $100812.10$ & $H_{25}$ & $448639.54$ & $H_{35}$ & $1209979.64$ & $H_{45}$ & $2546782.89$ \\
  $H_{6}$ & $7413.61$ & $H_{16}$ & $121604.52$ & $H_{26}$ & $503491.58$ & $H_{36}$ & $1315074.79$ & $H_{46}$ & $2718354.15$ \\
  $H_{7}$ & $11403.37$ & $H_{17}$ & $145073.90$ & $H_{27}$ & $562640.57$ & $H_{37}$ & $1426103.39$ & $H_{47}$ & $2897462.38$\\
  $H_{8}$ & $16612.10$ & $H_{18}$ & $171382.24$ & $H_{28}$ & $626248.53$ & $H_{38}$ & $1543210.96$ & $H_{48}$ & $3084269..56$ \\
  $H_{9}$ & $23201.79$ & $H_{19}$ & $200691.54$ & $H_{29}$ & $694477.44$ & $H_{39}$ & $1666559.50$ & $H_{49}$ & $3278937.71$ \\
  $H_{10}$ & $31334.44$ & $H_{20}$ & $233163.80$ & $H_{30}$ & $767489.32$ & $H_{40}$ & $1796310.99$ & $H_{50}$ & $3481628.82$ \\
  \hline
\end{tabular}
\end{table}

\begin{table}[htbp]
\setlength{\abovecaptionskip}{0.05cm}
 \centering \vspace{.3cm}
\caption{Complexity of linear heptagonal networks from $H_{1}$ to $H_{12}$.}
\begin{tabular}{cccccccc}
  \hline
  $H_n$ & $\tau(H_n)$ & $H_n$ & $\tau(H_n)$ & $H_n$ & $\tau(H_n)$ & $H_n$ & $\tau(H_n)$ \\
  \hline
  $H_{1}$ & $45$ & $H_{4}$ & $973080$ & $H_{7}$ & $21034100000$ & $H_{10}$ & $454673000000000$  \\
  $H_{2}$ & $1254$ & $H_{5}$ & $27106500$ & $H_{8}$ & $585934000000$ & $H_{11}$ & $12665600000000000$ \\
  $H_{3}$ & $34932$ & $H_{6}$ & $755090000$ & $H_{9}$ & $16322000000000$ & $H_{12}$ & $352817000000000000$ \\
  \hline
\end{tabular}
\end{table}

\section*{Acknowledgments}

This work is partially supported by National Natural Science
Foundation of China (Nos. 11601006, 11801007).


\end{document}